\documentclass[11pt]{amsart}

\usepackage[T1]{fontenc}
\usepackage[utf8]{inputenc}
\usepackage{amsmath,amssymb,amsthm}
\usepackage{cite}
\usepackage[margin=1in]{geometry}
\usepackage[hidelinks]{hyperref}

\numberwithin{equation}{section}

\newtheorem{theorem}{Theorem}
\newtheorem{conjecture}{Conjecture}
\theoremstyle{remark}
\newtheorem{remark}{Remark}

\begin{document}

\title[Belinsky's conjecture]{A Counterexample to Belinsky's Conjecture on
Ces\`aro Means at Lebesgue Points}
\author{Ushangi Goginava}
\address{Department of Mathematical Sciences, United Arab Emirates
University, P.O. Box 15551, Al Ain, Abu Dhabi, United Arab Emirates}
\email{zazagoginava@gmail.com}
\email{ugoginava@uaeu.ac.ae}
\date{}

\begin{abstract}
In 1997, Belinsky conjectured that, for convex subsequences, the logarithmic
growth condition of Carleson, Trigub, and Zagorodni\u{\i} is necessary and
sufficient for the arithmetic means of subsequential Fourier partial sums to
converge at every Lebesgue point of every integrable function. We disprove the
sufficiency part of this conjecture. More precisely, we construct a strictly
convex increasing sequence $(a_m)$ satisfying $a_m\leq 7m^8$ and a function
$f\in L^1(\mathbb T)$ for which $0$ is a Lebesgue point, $f(0)=0$, and the
means $m^{-1}\sum_{k=1}^m S_{a_k}f(0)$ are unbounded.
\end{abstract}

\maketitle

\section{Introduction}

Let $\mathbb T=[-\pi,\pi)$ , and let $S_nf(x)$
denote the $n$th partial sum of the Fourier series of a function
$f\in L^1(\mathbb T)$. Thus
\begin{equation*}
 S_nf(x)=\frac{1}{2\pi}\int_{-\pi}^{\pi}f(x+t)D_n(t)\,dt,
 \qquad
 D_n(t)=\frac{\sin((n+1/2)t)}{\sin(t/2)}.
\end{equation*}
A point $x\in\mathbb T$ is a Lebesgue point of $f$ if
\begin{equation*}
 \lim_{h\to0+}\frac{1}{2h}\int_{x-h}^{x+h}|f(t)-f(x)|\,dt=0.
\end{equation*}

In 1936, Zalcwasser \cite{Zalcwasser1936} asked how sparse a sequence of
positive integers $(a_k)$ may be while still preserving the convergence
\begin{equation}
 T_m^a f(x):=\frac{1}{m}\sum_{k=1}^m S_{a_k}f(x)\longrightarrow f(x),
 \qquad m\to\infty.
 \label{T}
\end{equation}
The study of this problem naturally separates into three directions:
\begin{enumerate}
 \item[(a)] uniform convergence in \eqref{T} for every continuous function;
 \item[(b)] almost-everywhere convergence in \eqref{T} for every integrable
 function;
 \item[(c)] convergence in \eqref{T} at every Lebesgue point of an integrable
 function.
\end{enumerate}

The first direction is completely solved for convex sequences. The uniform convergence of Ces\`aro means  along special subsequences  was invesigated by Bugrov
\cite{Bug},  Baiarstanova \cite{Baiar}, and
Long \cite{Long}. Later,  Carleson \cite{Carl}, and, independently, Zagorodni\u{\i} and
Trigub \cite{ZT} proved that, for a convex sequence $(a_k)$, the condition
\begin{equation}
 \sup_{n\geq1}\frac{\log a_n}{\sqrt n}<\infty
 \label{CTZ}
\end{equation}
is necessary and sufficient for uniform convergence in \eqref{T} for every
continuous function. This  questions is closely connected with
Salem's  work on strong summability \cite{Sal}.

The almost-everywhere problem in direction (b) also originates in the work of
Zalcwasser \cite{Zalcwasser1936}, who asked whether \eqref{T} holds almost
everywhere for every subsequence $(a_k)$ and every integrable function. G\'at
\cite{gat2023} answered this question in the negative by proving a general
almost-everywhere divergence result. He also identified a weak lacunarity
condition under which almost-everywhere convergence holds for every
integrable function \cite{gatConstr}; see also \cite{gatZ}. A complete
characterization of the subsequences having this universal almost-everywhere
convergence property remains open.

We now turn to direction (c). When $a_k=k$, convergence at every Lebesgue
point of every integrable function is the classical Fej\'er--Lebesgue theorem;
see \cite{Zygmund}. Related classical results on the summability of Fourier
series at Lebesgue points were obtained by Efimov \cite{Ef} and Faddeev
\cite{Fad}; see also the monograph of Weisz \cite{Weisz}.

Zalcwasser \cite{Zalcwasser1936} investigated \eqref{T} at Lebesgue points
for the quadratic subsequence $a_k=k^2$ and asked whether the same conclusion
holds for the more general power subsequences $a_k=k^r$. This question was
recently solved in \cite{Gog}.

Belinsky \cite{Belinsky1984,Belinsky1997} studied convergence at Lebesgue points along
non-polynomial subsequences. He constructed a sequence
$a_k\asymp e^{\sqrt[3]{k}}$ for which \eqref{T} converges at every Lebesgue
point of every integrable function. Motivated by this result and by the
Carleson--Trigub--Zagorodni\u{\i} theorem, he proposed the following
conjecture.

\begin{conjecture}[Belinsky, 1997]
Let $(a_k)$ be a convex increasing sequence of positive integers. Then
\eqref{CTZ} is necessary and sufficient for
\begin{equation*}
 T_m^a f(x)\longrightarrow f(x)
\end{equation*}
at every Lebesgue point $x$ of every function $f\in L^1(\mathbb T)$.
\end{conjecture}

The purpose of this paper is to disprove the sufficiency assertion in this
conjecture. Our main result is the following.

\begin{theorem}
\label{main}
There exist a strictly convex increasing sequence $(a_m)$ of positive
integers and a function $f\in L^1(\mathbb T)$ such that
\begin{enumerate}
 \item[(a)] $a_m\leq7m^8$ for every $m\geq1$;
 \item[(b)] $0$ is a Lebesgue point of $f$ and $f(0)=0$;
 \item[(c)]
 \begin{equation*}
  \sup_{m\geq1}\left|\frac{1}{m}\sum_{k=1}^mS_{a_k}f(0)\right|=\infty.
 \end{equation*}
\end{enumerate}
\end{theorem}

\section{Proof of Theorem \ref{main}}

\begin{proof}
\textbf{Step 1. Construction of the subsequence.}
We construct the sequence $(a_m)$ recursively. Let
\begin{equation*}
 a_1=1,\qquad a_2=4,\qquad a_3=9,\qquad a_4=16,
\end{equation*}
and call these terms block 1. Set $N_1=4$, the last index of block 1.
Suppose that the sequence has been constructed through the index $N_l$.
Define the next block by
\begin{equation}
 a_{N_l+s}=a_{N_l}+(4a_{N_l}+2)\frac{s(s+1)}{2},
 \qquad 1\leq s\leq N_l^2,
 \label{rec1}
\end{equation}
and set
\begin{equation*}
 N_{l+1}=N_l+N_l^2.
\end{equation*}
Thus block $l+1$ consists of the $N_l^2$ terms
\begin{equation*}
 a_{N_l+1},\ldots,a_{N_{l+1}}.
\end{equation*}
In particular, the beginning of the sequence has the block structure
\begin{equation*}
 \underbrace{a_1,\ldots,a_4}_{\text{block 1}},
 \quad
 \underbrace{a_5,\ldots,a_{4+4^2}}_{\text{block 2}},
 \quad\ldots.
\end{equation*}
Formula \eqref{rec1} also shows that every $a_m$ is a positive integer.

\medskip
\textbf{Step 2. Strict convexity of $(a_m)$.}
Set $d_m=a_{m+1}-a_m$. The initial differences satisfy
\begin{equation*}
 d_1=3<d_2=5<d_3=7<d_4=4a_4+2=66.
\end{equation*}
For $0\leq s<N_l^2$, formula \eqref{rec1} gives
\begin{align*}
 d_{N_l+s}
 &=a_{N_l+s+1}-a_{N_l+s}\\
 &=(4a_{N_l}+2)(s+1).
\end{align*}
Hence the differences increase strictly inside each block. At the boundary
between two consecutive blocks,
\begin{equation*}
 d_{N_{l+1}-1}=(4a_{N_l}+2)N_l^2,
\end{equation*}
whereas
\begin{align*}
 d_{N_{l+1}}
 &=4a_{N_{l+1}}+2\\
 &=(4a_{N_l}+2)\bigl(1+2N_l^2(N_l^2+1)\bigr)\\
 &>(4a_{N_l}+2)N_l^2
 =d_{N_{l+1}-1}.
\end{align*}
Therefore $d_{m+1}>d_m$ for every $m\geq1$, and $(a_m)$ is strictly
convex. In particular, it is strictly increasing.

\medskip
\textbf{Step 3. Polynomial growth of $(a_m)$.}
We first prove that
\begin{equation}
 a_{N_l}\leq N_l^6
 \label{ineq}
\end{equation}
for every $l\geq1$. For $l=1$, we have $a_{N_1}=16<4^6$. Suppose that
\eqref{ineq} holds for some $l$. Since $N_{l+1}=N_l(N_l+1)$, we obtain
\begin{align*}
 a_{N_{l+1}}
 &=a_{N_l}+(4a_{N_l}+2)
   \frac{N_l^2(N_l^2+1)}{2}\\
 &<a_{N_l}+3a_{N_l}N_l^2(N_l^2+1)\\
 &<a_{N_l}+6a_{N_l}N_l^4\\
 &\leq7N_l^{10}<N_l^{12}
 <\bigl(N_l(N_l+1)\bigr)^6=N_{l+1}^6.
\end{align*}
Thus \eqref{ineq} follows by induction.

Now let $m=N_l+s$, where $1\leq s\leq N_l^2$. Using \eqref{rec1} and
\eqref{ineq}, we find
\begin{align*}
 a_m
 &=a_{N_l}+(4a_{N_l}+2)\frac{s(s+1)}{2}\\
 &<a_{N_l}+3a_{N_l}s(s+1)\\
 &\leq7a_{N_l}s^2
 \leq7N_l^6s^2
 \leq7(N_l+s)^8=7m^8.
\end{align*}
This proves part (a) of the theorem.

\medskip
\textbf{Step 4. Construction of the function.}
For $v,j\geq1$, set
\begin{equation*}
 M_v=N_v+N_v^2,\qquad
 A_v=a_{M_v}+1,\qquad
 x_v=\frac{4\pi}{4a_{N_v}+2},\qquad
 c_j=2^{-j},\qquad
 R_j=2^{2j}.
\end{equation*}
Since $x_v\to0$, we may choose $v_1$ so that
\begin{equation*}
 \beta_1:=1,\qquad 2^{R_1}x_{v_1}<\beta_1.
\end{equation*}
Having chosen $v_j$ and $\beta_j$, define
\begin{equation}
 \beta_{j+1}
 =\min\left\{\frac{\beta_j}{4},x_{v_j},\frac{1}{A_{v_j}}\right\}
 \label{beta}
\end{equation}
and choose $v_{j+1}>v_j$ so large that
\begin{equation}
 2^{R_{j+1}}x_{v_{j+1}}<\beta_{j+1}.
 \label{seq2}
\end{equation}

For $0\leq r<R_j$, put
\begin{equation*}
 \rho_j=\frac{1}{16A_{v_j}},\qquad
 t_{j,r}=x_{v_j}\left(2^r+\frac12\right),\qquad
 H_{j,r}=[t_{j,r}-\rho_j,t_{j,r}+\rho_j],
\end{equation*}
and define
\begin{equation*}
 \sigma_{j,r}=
 \begin{cases}
  -1,&r=0,\\
  1,&1\leq r<R_j.
 \end{cases}
\end{equation*}
We shall use the following elementary support properties. Since
\begin{equation}
 \frac{\rho_j}{x_{v_j}}
 =\frac{4a_{N_{v_j}}+2}{64\pi(a_{M_{v_j}}+1)}
 <\frac{1}{16\pi}<\frac12,
 \label{rho-x}
\end{equation}
and $R_j\geq4$, conditions \eqref{seq2} and \eqref{rho-x} imply
\begin{equation}
 \bigcup_{r=0}^{R_j-1}H_{j,r}
 \subset(x_{v_j},\beta_j)\subset(0,\pi).
 \label{support}
\end{equation}

Define the $j$th block function by
\begin{equation*}
 u_j(t)=c_jA_{v_j}\sum_{r=0}^{R_j-1}
 \sigma_{j,r}t_{j,r}\mathbf 1_{H_{j,r}}(t).
\end{equation*}
As usual, $A\lesssim B$ means that $A\leq CB$ for a constant $C>0$
independent of the relevant parameters.
Using $|H_{j,r}|=2\rho_j=1/(8A_{v_j})$, we obtain
\begin{align}
 \|u_j\|_{L^1(0,\pi)}
 &\leq c_jA_{v_j}\sum_{r=0}^{R_j-1}t_{j,r}|H_{j,r}| \notag\\
 &=\frac{c_jx_{v_j}}{8}
   \left(2^{R_j}-1+\frac{R_j}{2}\right) \notag\\
 &\leq\frac{3}{16}c_j2^{R_j}x_{v_j}
 <\frac{3}{16}c_j\beta_j
 \lesssim2^{-j}.
 \label{uj}
\end{align}

For a sequence of signs $\varepsilon_j\in\{-1,1\}$ to be chosen in Step 6,
set
\begin{equation*}
 g(t)=\sum_{j=1}^{\infty}\varepsilon_ju_j(t),\qquad 0<t<\pi,
\end{equation*}
and define
\begin{equation*}
 f(0)=0,\qquad f(t)=g(|t|),\qquad -\pi<t<\pi,\quad t\neq0.
\end{equation*}
Estimate \eqref{uj} shows that the series defining $g$ converges in
$L^1(0,\pi)$. Moreover,
\begin{equation*}
 \|f\|_{L^1(\mathbb T)}
 =2\|g\|_{L^1(0,\pi)}
 \leq2\sum_{j=1}^{\infty}\|u_j\|_{L^1(0,\pi)}<\infty.
\end{equation*}

\medskip
\textbf{Step 5. The origin is a Lebesgue point of $f$.}
Since $\beta_{j+1}\leq\beta_j/4$, we have $\beta_j\to0$. For every
sufficiently small $h>0$, choose $j$ so that
\begin{equation}
 \beta_{j+1}<h\leq\beta_j.
 \label{h1}
\end{equation}
If $k<j$, then
\begin{equation*}
 h\leq\beta_j\leq\beta_{k+1}\leq x_{v_k},
\end{equation*}
and \eqref{support} gives $u_k=0$ on $[0,h]$. Therefore
\begin{align}
 \frac{1}{2h}\int_{-h}^{h}|f(t)-f(0)|\,dt
 &=\frac{1}{h}\int_0^h|g(t)|\,dt \notag\\
 &\leq\frac{1}{h}\int_0^h|u_j(t)|\,dt
 +\frac{1}{h}\sum_{k=j+1}^{\infty}\int_0^h|u_k(t)|\,dt.
 \label{lab1}
\end{align}

We first estimate the contribution of $u_j$. If
$H_{j,r}\cap[0,h]\neq\varnothing$, then \eqref{rho-x} and
$t_{j,r}\geq3x_{v_j}/2$ give
\begin{equation*}
 t_{j,r}\leq h+\rho_j<h+\frac23t_{j,r},
\end{equation*}
and hence $t_{j,r}<3h$. If there is at least one such interval, let
\begin{equation*}
 \overline r_j=\max\{r:t_{j,r}\leq3h\}.
\end{equation*}
Then $x_{v_j}2^{\overline r_j}\lesssim h$, and consequently
\begin{align*}
 \int_0^h|u_j(t)|\,dt
 &\lesssim c_jA_{v_j}
   \sum_{r:H_{j,r}\cap[0,h]\neq\varnothing}t_{j,r}|H_{j,r}|\\
 &\lesssim c_jx_{v_j}
   \sum_{r:t_{j,r}\leq3h}\left(2^r+\frac12\right)\\
 &\lesssim c_jx_{v_j}2^{\overline r_j}
 \lesssim hc_j.
\end{align*}
If no such interval exists, the same estimate is trivial. Thus
\begin{equation}
 \frac{1}{h}\int_0^h|u_j(t)|\,dt\lesssim c_j.
 \label{leb2}
\end{equation}

For the remaining blocks, \eqref{uj}, \eqref{beta}, and \eqref{h1} yield
\begin{align*}
 \sum_{k=j+1}^{\infty}\int_0^h|u_k(t)|\,dt
 &\leq\sum_{k=j+1}^{\infty}\|u_k\|_{L^1(0,\pi)}\\
 &\lesssim\sum_{k=j+1}^{\infty}c_k\beta_k\\
 &\lesssim c_j\beta_{j+1}
 \lesssim c_jh.
\end{align*}
Hence
\begin{equation}
 \frac{1}{h}\sum_{k=j+1}^{\infty}\int_0^h|u_k(t)|\,dt\lesssim c_j.
 \label{leb3}
\end{equation}
Combining \eqref{lab1}, \eqref{leb2}, and \eqref{leb3}, and observing that
$j\to\infty$ as $h\to0+$, we obtain
\begin{equation*}
 \frac{1}{2h}\int_{-h}^{h}|f(t)-f(0)|\,dt\lesssim c_j=2^{-j}\longrightarrow0.
\end{equation*}
Thus $0$ is a Lebesgue point of $f$.

\medskip
\textbf{Step 6. Divergence of the means.}
Let
\begin{equation*}
 K_m(t)=\frac{1}{m}\sum_{k=1}^mD_{a_k}(t),
 \qquad
 g_j(t)=\sum_{k=1}^j\varepsilon_ku_k(t).
\end{equation*}
Since $f$ and $K_m$ are even,
\begin{align}
 T_{M_{v_j}}^af(0)
 &=\frac{1}{2\pi}\int_{-\pi}^{\pi}f(t)K_{M_{v_j}}(t)\,dt \notag\\
 &=\frac{1}{\pi}\int_0^{\pi}g_j(t)K_{M_{v_j}}(t)\,dt
 +\frac{1}{\pi}\int_0^{\pi}(g(t)-g_j(t))K_{M_{v_j}}(t)\,dt.
 \label{T1}
\end{align}
The elementary estimate $\|D_n\|_{\infty}\leq2n+1$ gives
\begin{equation*}
 \|K_{M_{v_j}}\|_{\infty}
 \leq\frac{1}{M_{v_j}}\sum_{k=1}^{M_{v_j}}(2a_k+1)
 \leq2a_{M_{v_j}}+1<2A_{v_j}.
\end{equation*}
It follows from \eqref{uj} and \eqref{beta} that
\begin{align}
 \left|\int_0^{\pi}(g(t)-g_j(t))K_{M_{v_j}}(t)\,dt\right|
 &\lesssim A_{v_j}\sum_{k=j+1}^{\infty}\beta_kc_k \notag\\
 &\lesssim A_{v_j}c_{j+1}\beta_{j+1}
 \lesssim2^{-j}.
 \label{T2}
\end{align}

We now choose the signs $\varepsilon_j$. Suppose that
$\varepsilon_1,\ldots,\varepsilon_{j-1}$ have been chosen, and put
\begin{equation*}
 B_j=\int_0^{\pi}g_{j-1}(t)K_{M_{v_j}}(t)\,dt,
 \qquad g_0=0.
\end{equation*}
Choose
\begin{equation*}
 \varepsilon_j=
 \begin{cases}
  1,&B_j\geq0,\\
  -1,&B_j<0.
 \end{cases}
\end{equation*}
Then
\begin{equation}
 \left|\int_0^{\pi}g_j(t)K_{M_{v_j}}(t)\,dt\right|
 =\left||B_j|+\int_0^{\pi}u_j(t)K_{M_{v_j}}(t)\,dt\right|.
 \label{sg}
\end{equation}

Write
\begin{equation*}
 K_{M_{v_j}}(t)=\frac{Q_{M_{v_j}}(t)}{\sin(t/2)},
 \qquad
 Q_{M_{v_j}}(t)=\frac{1}{M_{v_j}}
 \sum_{k=1}^{M_{v_j}}\sin((a_k+1/2)t).
\end{equation*}
At the center $t_{j,r}$, split the last sum at $N_{v_j}$:
\begin{align}
 Q_{M_{v_j}}(t_{j,r})
 &=\frac{1}{M_{v_j}}\sum_{k=1}^{N_{v_j}}
   \sin((a_k+1/2)t_{j,r}) \notag\\
 &\quad+\frac{1}{M_{v_j}}\sum_{k=N_{v_j}+1}^{M_{v_j}}
   \sin((a_k+1/2)t_{j,r}).
 \label{I+II}
\end{align}
By \eqref{rec1},
\begin{align}
 &\sum_{k=N_{v_j}+1}^{M_{v_j}}
 \sin((a_k+1/2)t_{j,r}) \notag\\
 &\quad=\sum_{s=1}^{N_{v_j}^2}
 \sin\left(\pi\left(2^r+\frac12\right)
 +2s(s+1)\pi\left(2^r+\frac12\right)\right) \notag\\
 &\quad=N_{v_j}^2\sin\left(\pi\left(2^r+\frac12\right)\right)
 =N_{v_j}^2\cos(\pi2^r)
 =N_{v_j}^2\sigma_{j,r}.
 \label{I}
\end{align}
Here we used the fact that $s(s+1)$ is even. On the other hand,
\begin{equation}
 \sigma_{j,r}\sum_{k=1}^{N_{v_j}}
 \sin((a_k+1/2)t_{j,r})\geq-N_{v_j}.
 \label{II}
\end{equation}
Combining \eqref{I+II}--\eqref{II}, we obtain
\begin{equation}
 \sigma_{j,r}Q_{M_{v_j}}(t_{j,r})
 \geq\frac{N_{v_j}^2-N_{v_j}}{N_{v_j}^2+N_{v_j}}
 =\frac{N_{v_j}-1}{N_{v_j}+1}\geq\frac35.
 \label{s1}
\end{equation}

For $t\in H_{j,r}$, the mean value theorem gives
\begin{align}
 |Q_{M_{v_j}}(t)-Q_{M_{v_j}}(t_{j,r})|
 &\leq |t-t_{j,r}|\frac{1}{M_{v_j}}
   \sum_{k=1}^{M_{v_j}}(a_k+1/2)\\
 &\leq\rho_j(a_{M_{v_j}}+1/2)<\frac{1}{16}.
 \label{s2}
\end{align}
It follows from \eqref{s1} and \eqref{s2} that
\begin{equation}
 \sigma_{j,r}Q_{M_{v_j}}(t)
 >\frac35-\frac1{16}>\frac12,
 \qquad t\in H_{j,r}.
 \label{positive-Q}
\end{equation}

Using \eqref{support}, \eqref{positive-Q}, and
$\sin(t/2)\leq t/2$ for $0<t<\pi$, we conclude that
\begin{align*}
 \int_0^{\pi}u_j(t)K_{M_{v_j}}(t)\,dt
 &=c_jA_{v_j}\sum_{r=0}^{R_j-1}t_{j,r}
   \int_{H_{j,r}}\sigma_{j,r}
   \frac{Q_{M_{v_j}}(t)}{\sin(t/2)}\,dt\\
 &\geq c_jA_{v_j}\sum_{r=0}^{R_j-1}t_{j,r}
   \int_{H_{j,r}}\frac{dt}{t}\\
 &\geq c_jA_{v_j}\sum_{r=0}^{R_j-1}
   \frac{t_{j,r}}{t_{j,r}+\rho_j}|H_{j,r}|\\
 &\gtrsim c_jR_j=2^j.
\end{align*}
Consequently, \eqref{sg} yields
\begin{equation}
 \left|\int_0^{\pi}g_j(t)K_{M_{v_j}}(t)\,dt\right|\gtrsim2^j.
 \label{T3}
\end{equation}
Finally, combining \eqref{T1}, \eqref{T2}, and \eqref{T3}, we obtain
\begin{equation*}
 |T_{M_{v_j}}^af(0)|\gtrsim2^j-2^{-j}\longrightarrow\infty.
\end{equation*}
This proves part (c) and completes the proof.
\end{proof}

\begin{remark}
The preceding result leaves open the problem of characterizing those
increasing sequences $(a_k)$ for which \eqref{T} converges at every Lebesgue
point of every function in $L^1(\mathbb T)$.
\end{remark}

\section*{Declarations}

\noindent\textbf{Conflicts of interest.} The author declares that there are
no conflicts of interest.

\medskip \noindent\textbf{Data availability.} Not applicable.

\medskip \noindent\textbf{Funding.} This research received no external funding.

\end{document}